\documentclass[12pt]{amsart}
\usepackage[latin1]{inputenc}
\usepackage{amscd}
\usepackage{latexsym}
\usepackage{pstcol,pst-plot,pst-3d}
\usepackage{mathptmx}
\usepackage{multicol}
\psset{unit=0.7cm,linewidth=0.8pt,arrowsize=2.5pt 4}

\newpsstyle{fatline}{linewidth=1.5pt}
\newpsstyle{fyp}{fillstyle=solid,fillcolor=verylight}
\definecolor{verylight}{gray}{0.97}
\definecolor{light}{gray}{0.9}
\definecolor{medium}{gray}{0.85}
\def\frk{\frak}               

\def\Phi{{\frk n}}
\def\Phi{{\frk N}}
\def\opn#1#2{\def#1{\operatorname{#2}}} 
\opn\chara{char}
\opn\length{\ell}
\opn\pd{pd}
\opn\rk{rk}
\opn\projdim{proj\,dim}
\opn\injdim{inj\,dim}
\opn\rank{rank}
\opn\depth{depth}
\opn\grade{grade}
\opn\height{height}
\opn\embdim{emb\,dim}
\opn\codim{codim}

\opn\Tr{Tr}
\opn\bigrank{big\,rank}
\opn\superheight{superheight}\opn\lcm{lcm}
\opn\trdeg{tr\,deg}%
\opn\reg{reg}
\opn\lreg{lreg}
\opn\ini{in}
\opn\lpd{lpd}
\opn\div{div}
\opn\Div{Div}
\opn\cl{cl}
\opn\Cl{Cl}
\opn\Spec{Spec}
\opn\Supp{Supp}
\opn\supp{supp}
\opn\Sing{Sing}
\opn\Ass{Ass}
\opn\Ann{Ann}
\opn\Rad{Rad}
\opn\Soc{Soc}
\opn\Im{Im}
\opn\Ker{Ker}
\opn\Coker{Coker}
\opn\Am{Am}
\opn\Hom{Hom}
\opn\Tor{Tor}
\opn\Ext{Ext}
\opn\End{End}
\opn\Aut{Aut}
\opn\id{id}

\opn\nat{nat}
\opn\pff{pf}
\opn\Pf{Pf}
\opn\GL{GL}
\opn\SL{SL}
\opn\mod{mod}
\opn\ord{ord}
\opn\Gin{Gin}
\opn\aff{aff}
\opn\con{conv}
\opn\relint{relint}
\opn\st{st}
\opn\lk{lk}
\opn\cn{cn}
\opn\core{core}
\opn\vol{vol}
\opn\link{link}
\opn\star{star}
\opn\com{com}
\opn\gr{gr}

\def\pot#1#2{#1[\kern-0.28ex[#2]\kern-0.28ex]}

\opn\dirlim{\underrightarrow{\lim}}
\opn\inivlim{\underleftarrow{\lim}}
\let\sect=\cap

\let\iso=\cong
\let\Union=\bigcup
\let\Sect=\bigcap

\def\Implies{\ifmmode\Longrightarrow \else
       \unskip${}\Longrightarrow{}$\ignorespaces\fi}
\def\implies{\ifmmode\Rightarrow \else
       \unskip${}\Rightarrow{}$\ignorespaces\fi}
\def\iff{\ifmmode\Longleftrightarrow \else
       \unskip${}\Longleftrightarrow{}$\ignorespaces\fi}

\let\:=\colon
\newtheorem{Theorem}{Theorem}[section]
\newtheorem{Lemma}[Theorem]{Lemma}
\newtheorem{Corollary}[Theorem]{Corollary}

\let\epsilon\varepsilon
\let\phi=\varphi
\let\kappa=\varkappa
\def\qed{\ifhmode\textqed\fi
     \ifmmode\ifinner\quad\qedsymbol\else\dispqed\fi\fi}
\def\textqed{\unskip\nobreak\penalty50
      \hskip2em\hbox{}\nobreak\hfil\qedsymbol
      \parfillskip=0pt \finalhyphendemerits=0}
\def\dispqed{\rlap{\qquad\qedsymbol}}

\opn\dis{dis}
\def\pnt{{\raise0.5mm\hbox{\large\bf.}}}

\opn\Lex{Lex}

\begin{document}
\title{Cohen-Macaulay chordal graphs}
\author{J\"urgen Herzog, Takayuki Hibi and Xinxian Zheng}
\subjclass{13D02, 13P10, 13D40, 13A02}
\address{J\"urgen Herzog, Fachbereich Mathematik und
Informatik, Universit\"at Duisburg-Essen, Campus Essen,
45117 Essen, Germany}
\email{juergen.herzog@uni-essen.de}
\address{Takayuki Hibi, Department of Pure and Applied Mathematics,
Graduate School of Information Science and Technology,
Osaka University, Toyonaka, Osaka 560-0043, Japan}
\email{hibi@math.sci.osaka-u.ac.jp}
\address{Xinxian Zheng, Fachbereich Mathematik und
Informatik,
Universit\"at Duisburg-Essen, 45117 Essen, Germany}
\email{xinxian.zheng@uni-essen.de}
\date{}
\maketitle
\begin{abstract}
We classify all Cohen-Macaulay chordal graphs. In particular, it is shown that a chordal graph is Cohen-Macaulay if and only if it is unmixed.
\end{abstract}

\section*{Introduction}
To each finite graph $G$  with vertex set $[n]=\{1,\ldots, n\}$  and edge set $E(G)$ one associates the edge ideal $I(G)\subset K[x_1,\ldots, x_n]$ which is generated by all monomials $x_ix_j$ such that $\{i,j\}\in E(G)$. Here $K$ is an arbitrary field. The graph $G$ is called Cohen-Macaulay over $K$, if $K[x_1,\ldots, x_n]/I(G)$ is a Cohen-Macaulay ring, and is called Cohen-Macaulay if it is Cohen-Macaulay over any field.

Given a field $K$. The general problem is to classify the graphs which are Cohen-Macaulay over $K$. In this generality the problem is as hard as to classify all Cohen-Macaulay simplicial complexes,  because given a simplicial complex $\Delta$, one can
naturally construct a finite graph $G$ such that $G$ is
Cohen--Macaulay if and only if $\Delta$ is Cohen--Macaulay.  In fact,
if $P$ is the face poset of $\Delta$ (the poset consisting of all
faces of $\Delta$, ordered by inclusion), then $\Delta$ is
Cohen--Macaulay if and only if the order complex $\Delta(P)$ of $P$
is Cohen--Macaulay.  Since the order complex $\Delta(P)$ is flag,
i.e., every minimal nonface is a $2$-element subset, it follows that
there is a finite graph $G$ such that $I(G)$ coincides with the
Stanley--Reisner ideal of $\Delta(P)$.

Thus one cannot expect a general classification theorem. On the other hand, the first positive result was given by Villarreal \cite{V} who determined all Cohen-Macaulay trees. This result has been recently widely generalized in \cite{HH} where all bipartite Cohen-Macaulay graphs have been described. It turned out that the Cohen-Macaulay property of a bipartite graph does not depend on the field $K$.

In this note we classify all Cohen-Macaulay chordal graphs. Again it turns out that for chordal graphs the Cohen-Macaulay property is independent of the field $K$. Indeed we show that $G$ is Cohen-Macaulay if and only if the edge ideal $I(G)$ is height unmixed. One of our tools is Dirac's theorem \cite{D}  in a version as presented in \cite{HHZ}.

\section{Preliminaries}

Let $G$ be a finite graph on $[n]$ without loops, multiple edges
and isolated vertices, and $E(G)$ its edge set. The graph $G$ is called {\em chordal} if
all cycles of length $> 3$ has a chord.

A {\em stable
subset} or {\em clique} of $G$ is a subset $F$ of $[n]$ such that
$\{ i, j \} \in E(G)$ for all $i, j \in F$ with $i \neq j$. We
write $\Delta(G)$ for the simplicial complex on $[n]$ whose faces
are the stable subsets of $G$. For the proof of our main theorem we need the following property of chordal graphs \cite[Lemma 3.1]{HHZ} which  is related to  Dirac's theorem  \cite{D}.

\begin{Lemma}
\label{useful}
Let $G$ be a chordal graph. Then $\Delta(G)$ is a quasi-forest.
\end{Lemma}

We recall the definition of a quasi-forest introduced in \cite{Z}: let $\Delta$ be a simplicial complex, and $\mathcal{F}(\Delta)$ the set of its facets. A facet $F\in \mathcal{F}(\Delta)$ is called a {\em leaf}, if there exists a
facet $G$ (called a {\em branch} of $F$) with $G\neq F$ and such that $H\sect F\subset G\sect F$ for all $H\in{\mathcal F}(\Delta)$  with $H\neq F$.  We say that $\Delta$ is a {\em quasi-forest}, if there exists an order $F_1,\ldots, F_r$ of the facets of $\Delta$ such that for each $i=1,\ldots,r$, $F_i$ is a leaf of the simplicial complex $\langle F_1,\ldots, F_i\rangle$ (whose facets are $F_1,\ldots, F_i$).

\medskip

Let $K$ be a field. A graph $G$ is called {\em Cohen-Macaulay over $K$} if the edge ideal $I(G)=(\{x_ix_j \: \{i,j\}\in E(G)\})$ of $G$ is a Cohen-Macaulay ideal in $S=K[x_1,\ldots, x_n]$, in other words, if\ $S/I(G)$ is Cohen-Macaulay.

Suppose $G$ is Cohen-Macaulay over $K$. Then we say {\em $G$ is of type $r$} over $K$, if $r$ is the Cohen-Macaulay type of $S/I(G)$, that is, if $r$ is the minimal number of generators of the canonical module of $S/I(G)$.  The Cohen-Macaulay type of a Cohen-Macaulay ring $R$ can also be computed as the socle dimension of the residue class ring of $R$ modulo a  maximal regular sequence. The ring $R$ is Gorenstein, if and only if the Cohen-Macaulay type of $R$ is $1$. We say that {\em $G$ is Gorenstein over $K$}, if $S/I(G)$ is Gorenstein over $K$.

Finally we say that $G$  is {\em Cohen-Macaulay},   {\em of type $r$}, or {\em Gorenstein},  if $G$ has the corresponding  property over  any field.

The minimal prime ideals of $I(G)$ correspond to the minimal vertex covers of $G$. Recall that a {\em vertex cover} of $G$ is a subset $C\subset [n]$ such that $C\sect \{i,j\}\neq\emptyset$ for all $\{i,j\}\in E(G)$. It is called {\em minimal} if no proper subset of $C$ is a vertex cover of $G$. If we denote by ${\mathcal C}(G)$ the set of minimal vertex covers,  then the set of  ideals $\{(\{x_i\: i\in C\})\: C\in {\mathcal C}(G)\}$ is precisely the set of minimal prime ideals of $I(G)$.

Suppose again that $G$ is Cohen-Macaulay over $K$. Then the ideal $I(G)$ is height unmixed. Thus all minimal vertex covers of $G$   have the same cardinality.

\medskip
\noindent
For the proof of our main theorem we need the following algebraic fact:

\begin{Lemma}
\label{algebraic}
Let $R$ be a Noetherian ring, $S=R[x_1,\ldots, x_n]$ the polynomial ring over $R$, $k$ an integer with $0\leq k<n$, and $J$ the ideal $(I_1x_1,\ldots, I_kx_k, \{x_ix_j\}_{1\leq i<j\leq n})\subset S$, where $I_1,\ldots, I_k$ are ideals in $R$. Then the element $x=\sum_{i=1}^nx_i$ is a non-zerodivisor on $S/J$.
\end{Lemma}

\begin{proof}
For a subset $T\subset [n]$ we let $L_T$ be the ideal generated by all monomials $x_ix_j$ with $i,j\in T$ and $i<j$, and we set $I_T=\sum_{j\in T}I_j$ and $X_T=(\{x_j\}_{j\in T})$.

 It is easy to see that
\[
L_T=\Sect_{\ell\in T}X_{T\setminus\{\ell\}}.
\]
Hence we get
\begin{eqnarray*}
J&=& (I_1x_1,\ldots, I_kx_k, L_{[n]})
=\Sect_{T\subset [k]}(I_T,X_{[k]\setminus T}, L_{[n]})\\
&=& \Sect_{T\subset [k]}(I_T,X_{[k]\setminus T}, L_{[n]\setminus([k]\setminus T)})
= \Sect_{T\subset [k]\atop {\ell\in [n]\setminus([k]\setminus T)}}(I_T, X_{[k]\setminus T}, X_{([n]\setminus([k]\setminus T))\setminus \{\ell\}})\\
&=& \Sect_{T\subset [k]\atop {\ell\in [n]\setminus([k]\setminus T)}}(I_T, X_{[n]\setminus\{\ell\}}).
\end{eqnarray*}

Thus in order to prove  that $x$ is a non-zerodivisor modulo $J$ it suffices to show that $x$ is a non-zerodivisor modulo each of the ideals $(I_T, X_{[n]\setminus\{\ell\}})$. To see this we first pass to the residue class ring modulo $I_T$, and hence if we replace $R$ by $R/I_T$  it remains to be shown that $x$ is a non-zerodivisor on $R[x_1,\ldots, x_n]/(x_1,\ldots, x_{\ell-1},x_{\ell+1},\ldots x_n)$. But this is obviously the case.
\end{proof}

\section{The classification}

\begin{Theorem}
\label{main}

Let $K$ be a field, and let $G$ be a chordal graph on the vertex set $[n]$. Let $F_1,\ldots, F_m$ be the facets of $\Delta(G)$ which admit a free vertex. Then the following conditions are equivalent:
\begin{enumerate}
\item[(a)] $G$ is Cohen-Macaulay;
\item[(b)] $G$ is Cohen-Macaulay over $K$;
\item[(c)] $G$ is unmixed;
\item[(d)] $[n]$ is the disjoint union of $F_1,\ldots, F_m$.
\end{enumerate}
\end{Theorem}

\begin{proof}
(a)\implies (b) is trivial. 

(b)\implies (c): Since any Cohen-Macaulay ring is height unmixed it follows that $G$ is unmixed.

(c)\implies (d): 
Let $G$ be a unmixed chordal graph on $[n]$
and
$E(G)$ the set of edges of $G$.
Let $F_1, \ldots, F_m$  
denote the facets of $\Delta(G)$ with free vertices.
Fix a free vertex $v_i$ of $F_i$ 
and set $W = \{ v_1, \ldots, v_m \}$.
Suppose that $B = [n] \setminus (\Union_{i=1}^m F_i)
\neq \emptyset$
and write $G|_B$ for the induced subgraph of $G$ 
on $B$.
Since $\{ v_i, b \} \not\in E(G)$
for all $1 \leq i \leq m$ and for all $b \in B$,
if $X$ $(\subset B)$ is a minimal vertex cover 
of $G|_B$, then
$X \cup ((\Union_{i=1}^m F_i) \setminus W)$
is a minimal vertex cover of $G$.
In particular $G|_B$ is unmixed.
Since the induced subgraph $G|_B$ is
again chordal, by working with induction 
on the number of vertices, it follows that
if $H_1, \ldots, H_s$ are the facets of
$\Delta(G|_B)$ with free vertices, then
$B$ is the disjoint union 
$B = \Union_{j=1}^s H_j$.
Let $v'_j$ be a free vertex of $H_j$
and set $W' = \{ v'_1, \ldots, v'_s \}$.
Since
$((\Union_{i=1}^m F_i) \setminus W)
\cup (B \setminus W')$
is a minimal vertex cover of $G$ and
since $G$ is unmixed, 
every minimal vertex cover of $G$ consists of 
$n - (m + s)$ vertices.

We claim that $F_i \cap F_j = \emptyset$
if $i \neq j$.  
In fact, if, say, $F_1 \cap F_2 \neq \emptyset$
and if $w \in [n]$ satisfies
$w \in F_i$ for all $1 \leq i \leq \ell$, 
where $\ell \geq 2$,
and $w \not\in F_i$ for all $\ell < i \leq m$,
then
$Z = (\Union_{i=1}^m F_i) \setminus 
\{ w, v_{\ell + 1}, \ldots, v_m \}$
is a minimal vertex cover of
the induced subgraph $G' = G|_{[n] \setminus B}$
on $[n] \setminus B$.
Let $Y$ be a minimal vertex cover of $G$ with 
$Z \subset Y$.
Since $Y \cap B$ is a vertex cover of 
$G|_B$, one has $|Y \cap B| \geq |B| - s$.
Moreover, $|Y \cap ([n] \setminus B)| \geq
n - |B| - (m - \ell + 1) > n - |B| - m$.
Hence $|Y| > n - (m + s)$, a contradiction.

Consequently,
a subset $Y$ of $[n]$ is a minimal vertex cover 
of $G$ if and only if
$|Y \cap F_i| = |F_i| - 1$
for all $1 \leq i \leq m$ and
$|Y \cap H_j| = |H_j| - 1$
for all $1 \leq j \leq s$.

Now, since $\Delta(G|_B)$ is a quasi-forest,
one of the facets $H_1, \ldots, H_s$ must be a leaf
of $\Delta(G|_B)$.  Let, say, $H_1$ be a leaf
of $\Delta(G|_B)$.
Let $\delta$ and $\delta'$, where $\delta \neq \delta'$,
be free vertices of $H_1$
with $\{ \delta, a \} \in E(G)$
and $\{ \delta', a' \} \in E(G)$, where
$a$ and $a'$ belong to $[n] \setminus B$.
If $a \neq a'$ and if $\{ a, a' \} \in E(G)$, then
one has either
$\{ \delta, a' \} \in E(G)$
or $\{ \delta', a \} \in E(G)$,
because $G$ is chordal and 
$\{ \delta, \delta' \} \in E(G)$.
Hence there exists a subset 
$A \subset [n] \setminus B$ such that
\begin{enumerate}
\item[(i)] $\{ a, b \} \not\in E(G)$
for all $a, b \in A$ with $a \neq b$,
\item[(ii)] for each free vertex $\delta$
of $H_1$, one has $\{ \delta, a \} \in E(G)$
for some $a \in A$, and
\item[(iii)] for each $a \in A$, 
one has $\{ \delta, a \} \in E(G)$
for some free vertex $\delta$ of $H_1$.
\end{enumerate}
In fact, it is obvious that 
a subset $A \subset [n] \setminus B$ 
satisfying (ii) and (iii) exists.
If $\{ a, a' \} \in E(G)$,
$\{ \delta, a \} \in E(G)$ and
$\{ \delta, a' \} \not\in E(G)$ 
for some $a, a' \in A$ with $a \neq a'$
and for a free vertex $\delta$ of $H_1$,
then every free vertex $\delta'$ of $H_1$
with $\{ \delta', a' \} \in E(G)$
must satisfy $\{ \delta', a \} \in E(G)$.
Hence $A \setminus \{ a' \}$ satisfies
(ii) and (iii).
Repeating such the technique
yields a subset $A \subset [n] \setminus B$
satisfying (i), (ii) and (iii), as required. 

If $s>1$, then $H_1$ has a branch. Let $w_0 \not\in H_1$ be a vertex belonging to
a branch of the leaf $H_1$ of $\Delta(G|_B)$.
Thus $\{ \xi, w_0 \} \in E(G)$ for all nonfree
vertices $\xi$ of $H_1$.
We claim that either $\{ a, w_0 \} \not\in E(G)$
for all $a \in A$,
or one has $a \in A$
with $\{ a, \xi \} \in E(G)$
for every nonfree vertices $\xi$ of $H_1$.
To see why this is true,
if $\{ a, w_0 \} \in E(G)$ 
and $\{ \delta, a \} \in E(G)$
for some $a \in A$ and for some free vertex
$\delta$ of $H_1$, 
then one has a cycle 
$(a, \delta, \xi, w_0)$
of length four for every nonfree vertex $\xi$
of $H_1$.  Since $\{ \delta, w_0 \} \not\in E(G)$, 
one has $\{ a, \xi \} \in E(G)$.

Let $X$ be a minimal vertex cover of $G$ 
such that 
$X \subset [n] \setminus (A \cup \{ w_0 \})$
(resp. $X \subset [n] \setminus A$)
if $\{ a, w_0 \} \not\in E(G)$
for all $a \in A$
(resp. if one has $a \in A$
with $\{ a, \xi \} \in E(G)$
for every nonfree vertices $\xi$ of $H_1$.)
Then, for each vertex $\gamma$ of $H_1$,
there is $w \not\in X$ with $\{ \gamma, w \} \in E(G)$.
Hence $H_1 \subset X$, in contrast to our considerations before. 
This contradiction guarantees that $B = \emptyset$.
Hence $[n]$ is the disjoint union
$[n] = \Union_{i=1}^m F_i$, as required.

Finally suppose that $s=1$. Then $H_1$ is the only facet of $\Delta(G|_B)$.  Then $X=\Union_{i=1}^m(F_i\setminus v_i)$ is a minimal free vertex cover $G$ with $H_1\subset X$, a contradiction.

(d)\implies (c): 
Let $F_1, \ldots, F_m$ denote 
the facets of $\Delta(G)$ with free vertices
and, for each $1 \leq i \leq m$, 
write $F_i$ for the set of vertices of $F_i$.
Given a minimal vertex cover $X \subset [n]$ 
of $G$, one has
$|X \cap F_i| \geq |F_i| - 1$
for all $i$
since $F_i$ is a clique of $G$.
If, however, for some $i$, one has
$|X \cap F_i| = |F_i|$, i.e.,
$F_i \subset X$,
then
$X \setminus \{ v_i \}$
is a vertex cover of $G$
for any free vertex $v_i$
of $F_i$.  
This contradicts the fact that 
$X$ is a minimal vertex cover of $G$.
Thus
$|X \cap F_i| = |F_i| - 1$
for all $i$.
Since $[n]$ is the disjoint union
$[n] = \Union_{i=1}^m F_i$,
it follows that $|X| = n - m$
and $G$ is unmixed, as desired.

(c) and (d)\implies (a):  We  know that $G$ is unmixed. Moreover, if $v_i\in F_i$ is a free vertex,  then  $[n]\setminus \{v_1,\ldots, v_m\}$ is a minimal vertex cover of $G$. In particular it follows that $\dim S/I(G)=m$.

For $i=1,\ldots,m$, we set $y_i=\sum_{j\in F_i}x_j$. We will show that $y_1,\ldots, y_m$ is a regular sequence on $S/I(G)$. This then yields that $G$ is Cohen-Macaulay.

Let $F_i=\{i_1,\ldots, i_k\}$, and assume that $i_{\ell+1},\ldots, i_k$ are the free vertices of $F_i$. Let  $G'\subset G$ be the induced subgraph of $G$ on the vertex set $[n]\setminus \{i_1,\ldots, i_k\}$. Then $I(G)=(I(G'),J_1x_{i_1},\\
 J_2x_{i_2},\ldots, J_\ell x_{i_\ell}, J)$, where $J_j=(\{x_r\: \{r,i_j\}\in E(G)\})$ for $j=1,\ldots,\ell$, and where $J=(\{x_{i_r}x_{i_s}\: 1\leq r<s\leq k\})$.

Since $[n]$ is the disjoint union of $F_1,\ldots, F_m$ it follows that all generators of the ideal $(I(G'), y_1,\ldots, y_{i-1})$ belong to $K[\{x_i\}_{i\in [n]\setminus F_i}]$. Thus if we set
\[
R=K[\{x_i\}_{i\in [n]\setminus F_i}]/(I(G'), y_1,\ldots, y_{i-1}),
\]
then
\[
(S/I(G))/(y_1,\ldots,y_{i-1})(S/I(G))\iso R[x_{i_1},\ldots, x_{i_k}]/(I_1x_{i_1},\ldots, I_{\ell}x_{i_\ell}, \{x_{i_r}x_{i_s}\: 1\leq r<s\leq k\}),
\]
where for each $j$, the ideal $I_j$ is the image of $J_j$  under the residue class map onto $R$. Thus  Lemma \ref{algebraic} implies that $y_i$ is regular on $(S/I(G))/(y_1,\ldots, y_{i-1})(S/I(G))$.
\end{proof}

Let $G$ be an arbitrary graph on the vertex set $[n]$. An {\em independent set of $G$} is a set $S\subset [n]$ such that $\{i,j\}\not\in  E(G)$ for all $i,j\in S$. With this notion we can describe the type of a Cohen-Macaulay chordal graph.

\begin{Corollary}
\label{Gorenstein}
Let $G$ be a chordal graph, and let $F_1,\ldots, F_m$ be the facets of $\Delta(G)$ which have a free vertex. Let $i_j$ be a free vertex of $F_j$ for $j=1,\ldots, m$, and let $G'$ be the induced subgraph of $G$ on  the vertex set $[n]\setminus\{i_1,\ldots, i_m\}$.  Then 
\begin{enumerate}
\item[(a)] the type of $G$, is the number of maximal independent subsets of $G'$;
\item[(b)] $G$ is Gorenstein, if and only if $G$ is a disjoint union of edges.
\end{enumerate}
\end{Corollary}

\begin{proof} (a) Let $F\subset [n]$ and $S=K[x_1,\ldots, x_n]$. We note that if $J$ is the ideal generated by the set of monomials $\{x_ix_j\:  i,j\in F\;  \text{and}\;  i<j\}$, and $x=\sum_{i\in F}x_i$, then for any $i\in F$ one has that 
$$(S/J)/x(S/J)\iso  S_i/(\{x_j\: j\in F, j\neq i\})^2,$$ 
where $S_i=K[x_1,\ldots, x_{i-1}, x_{i+1},\ldots, x_n]$.

Thus if we factor by a maximal regular sequence as in the proof of Theorem \ref{main} we obtain a 0-dimensional  ring of the form 
\[
A=T/(P_1^2,\ldots, P_m^2,I(G'')).
\]
Here $P_j=(\{x_k\: k\in F_j,\; k\neq i_j\})$,   $G''$ is the subgraph of $G$ consisting of all edges which do not belong to any $F_j$, and  $T$ is the polynomial ring over $K$ in the set of variables $X=\{x_k\: k\in [n], k\neq i_j\; \text{for all}\; j=1,\ldots, m\}$. It is obvious that $A$  is obtained from the polynomial ring $T$ by factoring out the squares of all variables of $T$ and all $x_ix_j$ with $\{i,j\}\in E(G')$. Therefore $A$ has a $K$-basis of squarefree monomials corresponding to the independent subsets of $G'$, and  the socle  of $A$ is  generated as a $K$-vector space by  the monomials corresponding to the maximal independent subsets of $G'$.

(b) If $G$ is a disjoint union of edges, then $I(G)$ is a complete intersection, and hence Gorenstein.

Conversely, suppose that $G$ is Gorenstein. Then $A$ is Gorenstein. Since  $A$ a $0$-di\-men\-sio\-nal ring with monomial relations, $A$ is Gorenstein if and only if $A$ is a complete intersection. This is the case only if $E(G')=\emptyset$, in which case  $G$ is a disjoint union of egdes.
\end{proof}

\end{document}